\documentclass[12pt]{amsart}

\usepackage{amsmath}
\usepackage{amssymb}
\usepackage{amsthm}
\usepackage{color}
\usepackage{graphicx}
\usepackage{epsfig}

\newtheorem{theorem}{Theorem}[section]
\newtheorem{prop}[theorem]{Proposition}
\newtheorem{lemma}[theorem]{Lemma}
\newtheorem{cor}[theorem]{Corollary}

\theoremstyle{definition}

\newtheorem{rem}[theorem]{Remark}

\newcommand{\End}{\mathrm{End}}
 \newcommand{\Sh}{\mathrm{Sh}}
 \newcommand{\id}{{\mathit{id}}}
\newcommand{\meet}{\wedge}

\newcommand{\map}[1]{\xrightarrow{#1}}

\newcommand{\mat}[1]{\begin{pmatrix}#1\end{pmatrix}} 

\newcommand{\calA}{\mathcal{A}}
\newcommand{\calB}{\mathcal{B}}
\newcommand{\calK}{\mathcal{K}}
\newcommand{\calP}{\mathcal{P}}
\newcommand{\calS}{\mathcal{S}}
\newcommand{\calZ}{\mathcal{Z}}
\newcommand{\Sym}{{\mathit{Sym}}}
\newcommand{\SSym}{{\calS}\Sym}
\newcommand{\Z}{\mathbb{Z}}
\newcommand{\C}{\mathbb{C}}
\newcommand{\field}{\Bbbk} 
\newcommand{\kX}{\field\langle\!\langle X\rangle\!\rangle}

\newcommand{\red}[1]{\color{red}{#1}\color{black}}

\begin{document}
\title[Hopf algebra of uniform block permutations]{The Hopf algebra of uniform block permutations.\\
Extended abstract} 
\author{Marcelo Aguiar}
\address{Department of Mathematics\\
         Texas A\&M University\\
         College Station\\
         TX  77843\\
         USA}
\email{maguiar@math.tamu.edu}
\urladdr{http://www.math.tamu.edu/$\sim$maguiar}

\author{Rosa C. Orellana}
\address{Department of Mathematics\\
        Dartmouth College \\
        Hanover, NH 03755\\
        USA}
\email{Rosa.C.Orellana@Dartmouth.EDU}
\urladdr{http://www.math.dartmouth.edu/$\sim$orellana/}

\thanks{Aguiar supported in part by NSF grant DMS-0302423}
\thanks{Orellana supported in part by the Wilson Foundation}

\date{November 17, 2004}


 \maketitle

\begin{center}
 \begin{minipage}[c]{380pt}
 \small
 {\sc Abstract.}
 We introduce the Hopf algebra of uniform block permutations and show that it is  self-dual, free, and cofree. These results are closely related to
 the fact that  uniform block permutations form a factorizable inverse monoid.
This Hopf algebra
contains the Hopf algebra of permutations of Malvenuto and Reutenauer 
and the Hopf algebra of symmetric functions in non-commuting variables of
Gebhard, Rosas, and Sagan.
\end{minipage}\vspace{20pt}\\ 

\begin{minipage}[c]{380pt}
 \small
 {\sc R\'esum\'e.}
Nous pr\'esentons l'alg\`ebre de Hopf des permutations de blocs uniformes  est d\'emontrons qu'elle est auto duale, libre et colibre. Ces r\'esultats sont li\'es au fait que les permutations de blocs uniformes constituent un mono\"\i de inverse
factorisable.
Cette alg\`ebre de Hopf contient l'alg\`ebre de Hopf des permutations de Malvenuto et  Reutenauer et l'alg\`ebre de Hopf des fonctions 
sym\'etriques \`a variables non commutatives de Gebhard,  Rosas, et  Sagan.
\end{minipage}
\end{center}

\section{Uniform block permutations}

\subsection{Set partitions}\label{S:setpar}

Let $n$ be a non-negative integer and let $[n]:=\{1,2,\ldots, n\}$. A \emph{set partition} of $[n]$ is
a collection  of non-empty disjoint subsets of $[n]$, called {\em blocks}, whose union is $[n]$.
For example, $\calA=\big\{ \{2,5,7\}\{1,3\}\{6,8\}\{4\}\big\}$, is a set partition of $[8]$ with $4$ blocks.
  We often specify a set partition by listing the blocks from left to right so that the sequence formed by the minima of the blocks is increasing, and by listing the
  elements within each block in increasing order.
For instance,  the set partition above will be denoted 
$\calA=\{1,3\}\{2,5,7\}\{4\}\{6,8\}$.
We use $\calA\vdash [n]$ to indicate  that $\calA$ is a set partition of $[n]$.  


The {\em type} of a set partition $\calA$ of $[n]$ is the partition of $n$ formed by the sizes of the blocks of $\calA$. The symmetric group $S_n$ acts on the set of set partitions of $[n]$:
 given $\sigma\in S_n$ and $\calA\vdash [n]$, $\sigma(\calA)$ is the set partition whose blocks are $\sigma(A)$ for $A\in\calA$.
The orbit of $\calA$ consists of those set partitions of the same type as $\calA$.
The stabilizer of $\calA$ consists of those permutations that preserve the blocks, or that permute blocks of the same size. Therefore, the number 
of set partitions of type $1^{m_1}2^{m_2}\ldots n^{m_n}$ ($m_i$ blocks of size $i$) is
\begin{equation}\label{E:setpar}
 \frac{n!}{m_1! \cdots m_n! (1!)^{m_1}\cdots (n!)^{m_n}}\,.
 \end{equation}


%
%


\subsection{The monoid of uniform block permutations}

The monoid (and the monoid algebra) of uniform block permutations has been studied by FitzGerald \cite{f} and Kosuda \cite{k00,k01} in analogy to the partition algebra of Jones and Martin \cite{j,m}. 

A \emph{block permutation} of $[n]$ consists of  two set partitions $\calA$ and $\calB$ of $[n]$ with the same number of blocks and
 a bijection $f:\calA\to\calB$. For example, if $n=3$, $f(\{1,3\})=\{3\}$ and $f(\{2\})=\{1,2\}$ then $f$ 
is a block permutation.   
A block permutation is called \emph{uniform} if it maps each block of $\calA$ to a block of $\calB$ of the same 
cardinality.  For example, $f(\{1,3\})=\{1,2\}$, $f(\{2\})=\{3\}$
is uniform. Each permutation  may be viewed as a uniform block permutation for which all blocks have cardinality $1$.   
In this paper we only consider block permutations that are uniform. 

To specify a  uniform block permutation $f:\calA\to\calB$ we must choose
two set partitions $\calA$ and $\calB$ of the same type $1^{m_1}\ldots n^{m_n}$
and for each $i$ a bijection between the $m_i$ blocks of size $i$ of $\calA$
and those of $\calB$. We deduce from~\eqref{E:setpar} that the total number of uniform block permutations  of $[n]$ is 
\begin{equation}
u_n:=\sum_{1^{m_1}\ldots n^{m_n}\vdash n}  \left(\frac{n!}{(1!)^{m_1}\cdots (n!)^{m_n}}\right)^2\frac{1}{m_1!\cdots m_n!}
\end{equation}
where the sum runs over all partitions of $n$. Starting at $n=0$, the first values are
\[1,1,3,16,131,1496,22482,\ldots\]
This is sequence A023998 in~\cite{sl}. These numbers and generalizations
are studied in~\cite{sps}; in particular, the following recursion is given in~\cite[equation (11)]{sps}:
\[u_{n+1}=\sum_{k=0}^n\binom{n}{k}\binom{n+1}{k}u_k\,,\quad u_0=1\,.\]

\begin{figure}[!h]
\centering
\begin{picture}(370,70)(0,0)
\put(30,10){\circle*{3}}\put(30,60){\circle*{3}}\put(29,11){\line(5,3){80}}
\put(25,0){\small{1}}\put(25,67){\small{1}}
\put(50,10){\circle*{3}}\put(50,60){\circle*{3}}\put(30,10){\line(1,0){20}}
\put(45,0){\small{2}}\put(45,67){\small{2}}
\put(70,10){\circle*{3}}\put(70,60){\circle*{3}}\put(50,10){\line(2,1){100}}
\put(65,0){\small{3}}\put(65,67){\small{3}}
\put(90,10){\circle*{3}}\put(90,60){\circle*{3}}\put(70,10){\line(-4,5){40}}
\put(85,0){\small{4}}\put(85,67){\small{4}}
\put(110,10){\circle*{3}}\put(110,60){\circle*{3}}\put(70,10){\line(2,5){20}}
\put(105,0){\small{5}}\put(105,67){\small{5}}
\put(130,10){\circle*{3}}\put(130,60){\circle*{3}}\put(70,60){\line(1,0){20}}
\put(125,0){\small{6}}\put(125,67){\small{6}}
\put(150,10){\circle*{3}}\put(150,60){\circle*{3}}\qbezier(70,60)(50,67)(30,60)
\put(145,0){\small{7}}\put(145,67){\small{7}}
\put(170,10){\circle*{3}}\put(170,60){\circle*{3}}\put(90,10){\line(-4,5){40}}
\put(165,0){\small{8}}\put(165,67){\small{8}}
\put(110,10){\line(-4,5){40}}\put(110,10){\line(1,0){20}}\put(150,10){\line(2,5){20}}
\put(170,10){\line(-4,5){40}}
\put(200,30){$\sim$}
\put(230,10){\circle*{3}}\put(230,60){\circle*{3}}\put(229,11){\line(5,3){80}}
\put(225,0){\small{1}}\put(225,67){\small{1}}
\put(250,10){\circle*{3}}\put(250,60){\circle*{3}}\put(230,10){\line(1,0){20}}
\put(245,0){\small{2}}\put(245,67){\small{2}}
\qbezier(310,60)(330,67)(350,60)
\put(270,10){\circle*{3}}\put(270,60){\circle*{3}}\put(250,10){\line(2,1){100}}
\put(265,0){\small{3}}\put(265,67){\small{3}}
\put(290,10){\circle*{3}}\put(290,60){\circle*{3}}\put(270,10){\line(-4,5){40}}
\put(285,0){\small{4}}\put(285,67){\small{4}}
\put(310,10){\circle*{3}}\put(310,60){\circle*{3}}\qbezier(270,10)(290,3)(310,10)
\put(305,0){\small{5}}\put(305,67){\small{5}}
\put(330,10){\circle*{3}}\put(330,60){\circle*{3}}\put(270,60){\line(1,0){20}}
\put(325,0){\small{6}}\put(325,67){\small{6}}
\put(350,10){\circle*{3}}\put(350,60){\circle*{3}}\qbezier(270,60)(250,67)(230,60)
\put(345,0){\small{7}}\put(345,67){\small{7}}
\put(370,10){\circle*{3}}\put(370,60){\circle*{3}}\put(290,10){\line(-4,5){40}}
\put(365,0){\small{8}}\put(365,67){\small{8}}
\put(330,10){\line(-4,5){40}}\put(310,10){\line(1,0){20}}\put(350,10){\line(2,5){20}}
\put(370,10){\line(-4,5){40}}
\end{picture}
\caption{Two graphs representing the same uniform block permutation}\label{F:twographs} 

\end{figure}
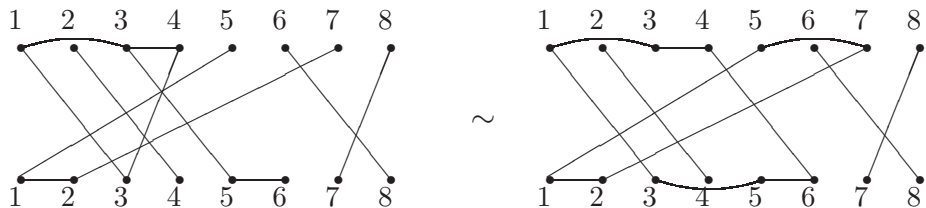

We represent uniform block permutations by means of  graphs.
For instance, either one of the two graphs in Figure~\ref{F:twographs}  represents the uniform block permutation $f$ given by
\[\{1,3,4\}\rightarrow\{3,5,6\},\ 
\{2\}\rightarrow \{4\}, \ \{5,7\}\rightarrow \{1,2\},\ \{6\}\rightarrow \{8\},\ \mbox{ and } 
\{8\}\rightarrow \{7\}\,.\]

Different graphs may represent the same uniform block permutation.
For a graph to represent a uniform block permutation $f:\calA\to\calB$ of $[n]$ the vertex set must consist
of two copies of $[n]$ (top and bottom) and each connected component must contain the same number of vertices on the top as on the bottom. The set partition $\calA$ is read off from the adjacencies on the top, $\calB$ from those on the bottom, and $f$ from those in between.

 The {\em diagram} of $f$ is the unique representing graph in
which all connected components are cycles and the elements in each cycle are joined in order,
as in the second graph of Figure~\ref{F:twographs}.

\begin{figure}[!htb]
\begin{picture}(350,70)(0,0)
\put(20,40){\large{$f =$}}
\put(50,60){\circle*{3}}\put(50,40){\circle*{3}}\put(50,40){\line(1,1){20}}
\put(70,60){\circle*{3}}\put(70,40){\circle*{3}}\put(70,40){\line(-1,1){20}}
\put(70,40){\line(1,0){20}}\put(90,40){\line(1,1){20}}\qbezier(50,60)(80,67)(110,60)
\put(90,60){\circle*{3}}\put(90,40){\circle*{3}}\put(110,40){\line(-1,1){20}}
\put(110,60){\circle*{3}}\put(110,40){\circle*{3}}\put(110,40){\line(1,0){20}}
\put(130,60){\circle*{3}}\put(130,40){\circle*{3}}\put(130,40){\line(1,1){20}}
\qbezier(90,60)(120,67)(150,60)
\put(150,60){\circle*{3}}\put(150,40){\circle*{3}}\put(150,40){\line(-1,1){20}}
\put(220,40){\large{$g =$}}
\put(250,60){\circle*{3}}\put(250,40){\circle*{3}}\put(250,40){\line(2,1){40}}
\put(270,60){\circle*{3}}\put(270,40){\circle*{3}}\put(270,40){\line(0,0){20}}
\put(290,60){\circle*{3}}\put(290,40){\circle*{3}}\put(290,40){\line(3,1){60}}
\put(310,60){\circle*{3}}\put(310,40){\circle*{3}}\qbezier(310,40)(330,30)(350,40)
\put(310,40){\line(-3,1){60}}\qbezier(250,60)(290,67)(330,60)
\put(330,60){\circle*{3}}\put(330,40){\circle*{3}}\put(330,40){\line(-1,1){20}}
\put(350,60){\circle*{3}}\put(350,40){\circle*{3}}\put(350,40){\line(-1,1){20}}
\end{picture}
\vskip 0in
\begin{picture}(350,80)(0,0)
\put(40,55){$g\cdot f=$}
\put(100,80){\circle*{3}}\put(100,60){\circle*{3}}\put(100,60){\line(1,1){20}}
\put(120,80){\circle*{3}}\put(120,60){\circle*{3}}\put(120,60){\line(-1,1){20}}
\put(120,60){\line(1,0){20}}\put(140,60){\line(1,1){20}}\qbezier(100,80)(130,87)(160,80)
\put(140,80){\circle*{3}}\put(140,60){\circle*{3}}\put(160,60){\line(-1,1){20}}
\put(160,80){\circle*{3}}\put(160,60){\circle*{3}}\put(160,60){\line(1,0){20}}
\put(180,80){\circle*{3}}\put(180,60){\circle*{3}}\put(180,60){\line(1,1){20}}
\qbezier(140,80)(170,87)(200,80)
\put(200,80){\circle*{3}}\put(200,60){\circle*{3}}\put(200,60){\line(-1,1){20}}
\put(100,45){\circle*{3}}\put(100,25){\circle*{3}}\put(100,25){\line(2,1){40}}
\put(100,45){\line(0,0){15}}
\put(120,45){\circle*{3}}\put(120,25){\circle*{3}}\put(120,25){\line(0,0){20}}
\put(120,45){\line(0,0){15}}
\put(140,45){\circle*{3}}\put(140,25){\circle*{3}}\put(140,25){\line(3,1){60}}
\put(140,45){\line(0,0){15}}
\put(160,45){\circle*{3}}\put(160,25){\circle*{3}}\qbezier(160,25)(180,15)(200,25)
\put(160,45){\line(0,0){15}}
\put(160,25){\line(-3,1){60}}\qbezier(100,45)(130,52)(180,45)
\put(180,45){\circle*{3}}\put(180,25){\circle*{3}}\put(180,25){\line(-1,1){20}}
\put(180,45){\line(0,0){15}}
\put(200,45){\circle*{3}}\put(200,25){\circle*{3}}\put(200,25){\line(-1,1){20}}
\put(200,45){\line(0,0){15}}
\put(220,55){$=$}

\put(250,65){\circle*{3}}\put(250,45){\circle*{3}}\put(250,45){\line(0,0){20}}
\put(250,45){\line(1,0){20}}\qbezier(250,65)(280,72)(310,65)
\put(270,65){\circle*{3}}\put(270,45){\circle*{3}}\put(270,45){\line(2,1){40}}
\put(290,65){\circle*{3}}\put(290,45){\circle*{3}}\put(290,45){\line(2,1){40}}
\put(310,65){\circle*{3}}\put(310,45){\circle*{3}}\put(310,45){\line(1,0){20}}
\put(310,45){\line(-2,1){40}}\put(270,65){\line(1,0){20}}\qbezier(290,65)(320,72)(350,65)
\put(330,65){\circle*{3}}\put(330,45){\circle*{3}}\put(330,45){\line(1,0){20}}
\put(350,65){\circle*{3}}\put(350,45){\circle*{3}}\put(350,45){\line(0,0){20}}
\end{picture}
\caption{Product of uniform block permutations}\label{F:product}
\end{figure}
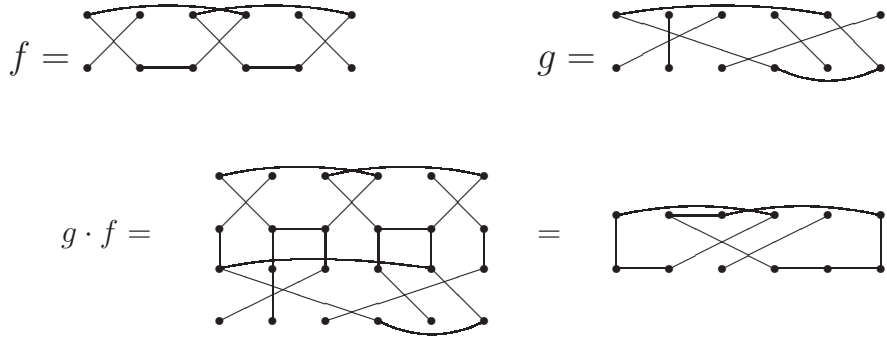

The set $P_n$ of block permutations of $[n]$ is a monoid.
The product $g\cdot f$ of   two uniform block permutations $f$ and $g$ of $[n]$ is
obtained by gluing the bottom of a graph representing $f$ to the top of a graph representing $g$. The resulting graph represents a uniform block permutation which does not depend on the graphs chosen.   
An example is given in Figure~\ref{F:product}. Note that gluing the diagram of $f$ to the diagram of $g$ may not result in the diagram of $g\cdot f$.

 The identity is the uniform
block permutation that maps $\{i\}$ to $\{i\}$ for all $i$. Viewing permutations as uniform block permutations as above, we get that the symmetric group $S_n$ is a submonoid of $P_n$.

We recall a presentation of the monoid $P_n$ given in \cite{f,k00,k01}.
Consider the uniform block permutations $b_i$ and $s_i$ with diagrams

\begin{picture}(375,80)(0,0)
\put(20,35){\large{$b_i =$}}
\put(50,60){\circle*{3}}\put(50,20){\circle*{3}}\put(50,20){\line(0,0){40}}
\put(45,67){\tiny{1}}\put(45,10){\tiny{1}}
\put(60,35){\Large{$\cdots$}}
\put(90,60){\circle*{3}}\put(90,20){\circle*{3}}\put(90,20){\line(0,0){40}}
\put(80,67){\tiny{$i-1$}}\put(80,10){\tiny{$i-1$}}
\put(105,60){\circle*{3}}\put(105,20){\circle*{3}}\put(105,20){\line(0,0){40}}
\put(103,67){\tiny{$i$}}\put(103,10){\tiny{$i$}}
\put(105,60){\line(1,0){15}}\put(105,20){\line(1,0){15}}
\put(120,60){\circle*{3}}\put(120,20){\circle*{3}}\put(120,20){\line(0,0){40}}
\put(115,67){\tiny{$i+1$}}\put(115,10){\tiny{$i+1$}}
\put(135,60){\circle*{3}}\put(135,20){\circle*{3}}\put(135,20){\line(0,0){40}}
\put(135,67){\tiny{$i+2$}}\put(135,10){\tiny{$i+2$}}
\put(145,35){\Large{$\cdots$}}
\put(175,60){\circle*{3}}\put(175,20){\circle*{3}}\put(175,20){\line(0,0){40}}
\put(173,67){\tiny{$n$}}\put(173,10){\tiny{$n$}}
\put(220,35){\large{$s_i =$}}
\put(250,60){\circle*{3}}\put(250,20){\circle*{3}}\put(250,20){\line(0,0){40}}
\put(245,67){\tiny{1}}\put(245,10){\tiny{1}}
\put(260,35){\Large{$\cdots$}}
\put(290,60){\circle*{3}}\put(290,20){\circle*{3}}\put(290,20){\line(0,0){40}}
\put(280,67){\tiny{$i-1$}}\put(280,10){\tiny{$i-1$}}
\put(305,60){\circle*{3}}\put(305,20){\circle*{3}}\put(304,21){\line(2,5){15}}
\put(303,67){\tiny{$i$}}\put(303,10){\tiny{$i$}}
\put(320,60){\circle*{3}}\put(320,20){\circle*{3}}\put(321,21){\line(-2,5){15}}
\put(315,67){\tiny{$i+1$}}\put(315,10){\tiny{$i+1$}}
\put(335,60){\circle*{3}}\put(335,20){\circle*{3}}\put(335,20){\line(0,0){40}}
\put(335,67){\tiny{$i+2$}}\put(335,10){\tiny{$i+2$}}
\put(345,35){\Large{$\cdots$}}
\put(375,60){\circle*{3}}\put(375,20){\circle*{3}}\put(375,20){\line(0,0){40}}
\put(373,67){\tiny{$n$}}\put(373,10){\tiny{$n$}}
\end{picture}

The monoid $P_n$ is generated by the 
elements $\{ b_i, s_i\, |\, 1\leq i\leq n-1\}$ subject to the following relations:
\begin{enumerate}
\item $s_i^2 =1, \qquad b_i^2 =b_i, \qquad 1\leq i \leq n-1$;
\item $s_is_{i+1}s_i=s_{i+1}s_is_{i+1}, \qquad s_ib_{i+1}s_i=s_{i+1}b_is_{i+1}, \qquad 1\leq i\leq n-2$;
\item $s_is_j=s_js_i, \qquad b_is_j=s_jb_i, \qquad |i-j|>1$;
\item $ b_is_i=s_ib_i=b_i, \qquad 1\leq i \leq n-1$;
\item $ b_ib_j=b_jb_i,\qquad 1\leq i,j\leq n-1$.
\end{enumerate}

The submonoid generated by the elements $s_i$, $1\leq i\leq n-1$ is the symmetric group $S_n$, viewed as a submonoid of $P_n$ as above.

We will see in Sections~\ref{S:selfdual} and~\ref{S:weak} that $P_n$ is a factorizable inverse monoid. Therefore, a presentation for $P_n$ may also be derived from
the results of~\cite{eef}.


\subsection{An ideal indexed by set partitions}\label{S:ideal}

Let $\field P_n$ be the monoid algebra of $P_n$ over a commutative ring $\field$.

Given a set partition $\calA\vdash [n]$,  let $Z_\calA\in\field P_n$ denote the sum of all uniform block permutations $f:\calA\to\calB$, where $\calB$ varies:
$$Z_\calA:=\sum_{f:\calA\to\calB}f\,.$$
For instance,
\vskip 0in 
\begin{picture}(400,45)(0,0)
\put(0,15){$Z_{\{1,3\}\{2,4\}}=$}
\put(70,10){\circle*{3}}\put(80,10){\circle*{3}}\put(90,10){\circle*{3}} \put(100,10){\circle*{3}}
\put(70,30){\circle*{3}}\put(80,30){\circle*{3}}\put(90,30){\circle*{3}} \put(100,30){\circle*{3}}
\qbezier(70,30)(80,37)(90,30)\qbezier(80,30)(90,37)(100,30)\put(70,10){\line(1,0){10}}\put(90,10){\line(1,0){10}}
\put(70,10){\line(0,0){20}}\put(80,10){\line(1,2){10}}\put(90,10){\line(-1,2){10}}\put(100,10){\line(0,0){20}}
\put(110,17){$+$}
\put(130,10){\circle*{3}}\put(140,10){\circle*{3}}\put(150,10){\circle*{3}} \put(160,10){\circle*{3}}
\put(130,30){\circle*{3}}\put(140,30){\circle*{3}}\put(150,30){\circle*{3}} \put(160,30){\circle*{3}}
\qbezier(130,30)(140,37)(150,30)\qbezier(140,30)(150,37)(160,30)\put(130,10){\line(0,0){20}}
\qbezier(130,10)(140,17)(150,10)\qbezier(140,10)(150,17)(160,10)\put(140,10){\line(0,0){20}}
\put(150,10){\line(0,0){20}}\put(160,10){\line(0,0){20}}
\put(170,17){$+$}
\put(190,10){\circle*{3}}\put(200,10){\circle*{3}}\put(210,10){\circle*{3}} \put(220,10){\circle*{3}}
\put(190,30){\circle*{3}}\put(200,30){\circle*{3}}\put(210,30){\circle*{3}} \put(220,30){\circle*{3}}
\qbezier(190,30)(200,37)(210,30)\qbezier(200,30)(210,37)(220,30)\put(190,10){\line(0,0){20}}
\qbezier(190,10)(205,17)(220,10)\put(200,10){\line(1,0){10}}\put(200,10){\line(0,0){20}}
\put(210,10){\line(1,2){10}}\put(220,10){\line(-1,2){10}}
\put(230,15){$+$}
\put(250,10){\circle*{3}}\put(260,10){\circle*{3}}\put(270,10){\circle*{3}} \put(280,10){\circle*{3}}
\put(250,30){\circle*{3}}\put(260,30){\circle*{3}}\put(270,30){\circle*{3}} \put(280,30){\circle*{3}}
\qbezier(250,30)(260,37)(270,30)\qbezier(260,30)(270,37)(280,30)\put(250,10){\line(1,2){10}}
\qbezier(250,10)(265,17)(280,10)\put(260,10){\line(1,0){10}}\put(260,10){\line(-1,2){10}}
\put(270,10){\line(0,0){20}}\put(280,10){\line(0,0){20}}
\put(290,15){$+$}
\put(310,10){\circle*{3}}\put(320,10){\circle*{3}}\put(330,10){\circle*{3}} \put(340,10){\circle*{3}}
\put(310,30){\circle*{3}}\put(320,30){\circle*{3}}\put(330,30){\circle*{3}} \put(340,30){\circle*{3}}
\qbezier(310,30)(320,37)(330,30)\qbezier(320,30)(330,37)(340,30)\put(310,10){\line(1,2){10}}
\qbezier(310,10)(320,17)(330,10)\qbezier(320,10)(330,17)(340,10)\put(320,10){\line(-1,2){10}}
\put(330,10){\line(1,2){10}}\put(340,10){\line(-1,2){10}}
\put(350,15){$+$}
\put(370,10){\circle*{3}}\put(380,10){\circle*{3}}\put(390,10){\circle*{3}} \put(400,10){\circle*{3}}
\put(370,30){\circle*{3}}\put(380,30){\circle*{3}}\put(390,30){\circle*{3}} \put(400,30){\circle*{3}}
\qbezier(370,30)(380,37)(390,30)\qbezier(380,30)(390,37)(400,30)\put(370,10){\line(1,2){10}}
\put(370,10){\line(1,0){10}}\put(390,10){\line(1,0){10}}\put(380,10){\line(1,1){20}}
\put(390,10){\line(-1,1){20}}\put(400,10){\line(-1,2){10}}
\end{picture}

\begin{lemma} Let $\calA$ be a set partition of $[n]$ and $\sigma$ a permutation of $[n]$. Then
 \[\sigma \cdot Z_\calA= Z_\calA \text{ \ and \ }Z_\calA \cdot \sigma = Z_{\sigma^{-1}(\calA)}\,.\]
 In addition,
\[Z_\calA\cdot b_i= \left\{ \begin{array}{ll} Z_\calA & \mbox{if $i$ and $i+1$ belong to the same block of $\calA$}\\
\rule{0pt}{20pt}{|A|+|A'|\choose |A|} Z_\calB  & \mbox{if $i$ and $i+1$ belong to different blocks $A$ and $A'$ of $\calA$}\end{array} \right. \] 
where the set partition $\calB$ is obtained by merging the blocks $A$ and $A'$ of $\calA$ and keeping the others. 
\end{lemma}


Let $\calZ_n$ denote the subspace of $\field P_n$ linearly spanned by the elements $Z_\calA$ as $\calA$ runs over all set partitions of $[n]$.

\begin{cor}\label{C:ideal}
$\calZ_n$ is a right ideal of the monoid algebra $\field P_n$. 
\end{cor}


\section{The Hopf algebra of uniform block permutations}
In this section we define the Hopf algebra of uniform block permutations. It contains the
Hopf algebra of permutations of Malvenuto and Reutenauer as a Hopf subalgebra.

\subsection{Schur-Weyl duality for uniform block permutations}

Let $r$ and $m$ be positive integers. Consider the complex reflection group
\[G(r,1,m):=\Z_r\wr S_m\,.\]
Let $t$ denote the generator of the cyclic group $\Z_r$, $t^r=1$.

Let $V$ be the {\em monomial} representation of $G(r,1,m)$. Thus, $V$ is
an $m$-dimensional vector space with a  basis $\{ e_1,e_2,\ldots, e_m\}$ on which  $G(r,1,m)$ acts as follows:
\[t\cdot e_1 = e^{2\pi i/r} e_1\,,\quad t\cdot e_i =e_i \text{ for  $i>1$, and \ }
\sigma\cdot e_i = e_{\sigma(i)} Ý\text{ for $\sigma \in S_m$.}\]

Consider now the diagonal action of $G(r,1,m)$ on the tensor powers $V^{\otimes n}$,
\[Ýg\cdot (e_{i_1}e_{i_2}\cdots e_{i_n}) =(g\cdot e_{i_1})(g\cdot e_{i_2}) \cdots (g\cdot e_{i_n})\,.\]
The centralizer of this representation has been calculated by Tanabe.
\begin{prop}~\cite{t} \label{P:block-duality}
There is a right action of the monoid $P_n$ on $V^{\otimes n}$ determined by
\[(e_{i_1}\cdots e_{i_n})\cdot b_j=\delta(i_j,i_{j+1})e_{i_1}\cdots e_{i_n} \text{ \ and \ } (e_{i_1}\cdots e_{i_n})\cdot\sigma= e_{i_{\sigma(1)}}\cdots e_{i_{\sigma(n)}}\]
for $1\leq i\leq n-1$ and   $\sigma\in S_n$. This action commutes with the left action of $G(r,1,m)$ on $V^{\otimes n}$.
Moreover, if $m\geq 2n$ and $r>n$ then the resulting map
\begin{equation} \label{E:block-duality}
\C P_n\to \End_{G(r,1,m)}(V^{\otimes n})
\end{equation}
is an isomorphism of algebras.
\end{prop}
Classical Schur-Weyl duality states that the symmetric group algebra can be similarly recovered from the diagonal action of $GL(V)$ on $V^{\otimes n}$: if $\dim V\geq n$ then
\begin{equation} \label{E:classical-duality}
\C S_n\cong \End_{GL(V)}(V^{\otimes n})\,.
\end{equation}
Malvenuto and Reutenauer~\cite{mr} deduce from here
 the existence of a multiplication among permutations as follows. 
 Given $\sigma\in S_p$ and $\tau\in S_q$, view them as linear endomorphisms
 of the tensor algebra
 \[T(V):=\bigoplus_{n\geq 0}V^{\otimes n}\]
 by means of~\eqref{E:classical-duality} ($\sigma$ acts as $0$ on $V^{\otimes n}$ if $n\neq p$, similarly for $\tau$). The tensor algebra is a Hopf algebra, so we can form the {\em convolution} product of any two linear endomorphisms:
 \[T(V)\map{\Delta}T(V)\otimes T(V)\map{\sigma\otimes\tau}T(V)\otimes T(V)\map{m}T(V)\,,\]
 where $\Delta$ and $m$ are the coproduct and product of the tensor algebra. 
 Since these two maps commute with the action of $GL(V)$, the convolution  of $\sigma$ and $\tau$ belongs to $\End_{GL(V)}(V^{\otimes n})$, where $n=p+q$.
 Therefore, there exist an element $\sigma\ast\tau\in\C S_n$ whose right action
 equals the convolution  of $\sigma$ and $\tau$. This is the product of Malvenuto and Reutenauer.
 
The same argument applies to uniform block permutations, in view of Proposition~\ref{P:block-duality}. We proceed to describe the resulting operation
in explicit terms. As for permutations, this structure can be enlarged to that
of a graded Hopf algebra.
 
\subsection{Product and coproduct of uniform block permutations}

Consider the graded vector space 
$$\calP:=\bigoplus_{n\geq 0} \field P_n\,.$$
$P_0$ consists of the unique uniform block permutation of $[n]$, represented by the empty diagram, which we denote by $\emptyset$.

\begin{figure}[!htb]
\begin{picture}(290,100)(0,0)
\put(40,75){$f=$} \put(70,90){\circle*{3}}\put(70,70){\circle*{3}}\put(70,70){\line(2,1){40}}
\put(90,90){\circle*{3}}\put(90,70){\circle*{3}}\put(90,70){\line(-1,1){20}}
\put(110,90){\circle*{3}}\put(110,70){\circle*{3}}\put(110,70){\line(-1,1){20}}
\put(130,90){\circle*{3}}\put(130,70){\circle*{3}}\put(110,70){\line(1,0){20}}
\put(130,70){\line(0,0){20}}\qbezier(90,90)(110,97)(130,90)
\put(68,95){\tiny{1}}\put(88,95){\tiny{2}}\put(108,95){\tiny{3}}\put(128,95){\tiny{4}}
\put(68,63){\tiny{1}}\put(88,63){\tiny{2}}\put(108,63){\tiny{3}}\put(128,63){\tiny{4}}
\put(135,70){\Large{,}}
\put(170,75){$g=$}\put(200,90){\circle*{3}}\put(200,70){\circle*{3}}\qbezier(200,70)(220,60)(240,70)
\put(240,70){\line(1,0){20}}\put(260,70){\line(1,1){20}}\put(220,90){\line(1,0){20}}
\put(220,90){\circle*{3}}\put(220,70){\circle*{3}}\put(200,70){\line(1,1){20}}
\put(240,90){\circle*{3}}\put(240,70){\circle*{3}}\qbezier(240,90)(260,97)(280,90)
\put(260,90){\circle*{3}}\put(260,70){\circle*{3}}\put(280,70){\line(-4,1){80}}
\put(280,90){\circle*{3}}\put(280,70){\circle*{3}}\put(220,70){\line(2,1){40}}
\put(198,95){\tiny{1}}\put(218,95){\tiny{2}}\put(238,95){\tiny{3}}\put(258,95){\tiny{4}}
\put(278,95){\tiny{5}}\put(198,63){\tiny{1}}\put(218,63){\tiny{2}}\put(238,63){\tiny{3}}
\put(258,63){\tiny{4}}\put(278,63){\tiny{5}}
\put(80,25){$f\times g=$}
\put(130,40){\circle*{3}}\put(130,20){\circle*{3}}\put(130,20){\line(2,1){40}}
\put(150,40){\circle*{3}}\put(150,20){\circle*{3}}\put(150,20){\line(-1,1){20}}
\put(170,40){\circle*{3}}\put(170,20){\circle*{3}}\put(170,20){\line(-1,1){20}}
\put(190,40){\circle*{3}}\put(190,20){\circle*{3}}\put(170,20){\line(1,0){20}}
\put(190,20){\line(0,0){20}}\qbezier(150,40)(170,47)(190,40)
\put(128,45){\tiny{1}}\put(148,45){\tiny{2}}\put(168,45){\tiny{3}}\put(188,45){\tiny{4}}
\put(208,45){\tiny{5}}\put(228,45){\tiny{6}}\put(248,45){\tiny{7}}\put(268,45){\tiny{8}}
\put(288,45){\tiny{9}}
\put(128,13){\tiny{1}}\put(148,13){\tiny{2}}\put(168,13){\tiny{3}}\put(188,13){\tiny{4}}
\put(208,13){\tiny{5}}\put(228,13){\tiny{6}}\put(248,13){\tiny{7}}\put(268,13){\tiny{8}}
\put(288,13){\tiny{9}}
\put(210,40){\circle*{3}}\put(210,20){\circle*{3}}\qbezier(210,20)(230,10)(250,20)
\put(250,20){\line(1,0){20}}\put(270,20){\line(1,1){20}}\put(230,40){\line(1,0){20}}
\put(230,40){\circle*{3}}\put(230,20){\circle*{3}}\put(210,20){\line(1,1){20}}
\put(250,40){\circle*{3}}\put(250,20){\circle*{3}}\qbezier(250,40)(270,47)(290,40)
\put(270,40){\circle*{3}}\put(270,20){\circle*{3}}\put(290,20){\line(-4,1){80}}
\put(290,40){\circle*{3}}\put(290,20){\circle*{3}}\put(230,20){\line(2,1){40}}
\end{picture}
\caption{Concatenation of diagrams} \label{F:concatenation}
\end{figure}
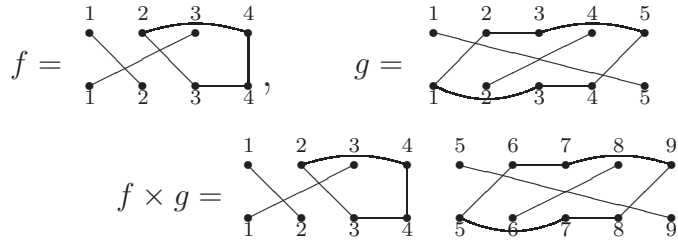

 Let $f$ and  $g$ be uniform block permutations of $[n]$ and $[m]$ respectively.
 Adding $n$ to every entry in the diagram of $g$ and placing it to the right of the
 diagram of $f$ we obtain the diagram of a uniform block permutation of $[n+m]$, called the concatenation of $f$ and $g$ and denoted  $f\times g$. Figure~\ref{F:concatenation} shows an example.

Let $\Sh(n,m)$ denote the set of  $(n,m)$-shuffles, that is, those 
 permutations $\xi \in S_{n+m}$ such that 
$$\xi(1)<\xi(2)< \cdots < \xi(n) \mbox{ \ and \ } \xi(n+1)<\xi(n+2)<\cdots < \xi(n+m)\,.$$
Let $sh_{n,m}\in \field S_{n+m}$ denote the sum of all $(n,m)$-shuffles.  

The product $\ast$ on 
$\calP$ is defined by 
$$f\ast g := sh_{n,m}\cdot (f\times g)\in \field P_{n+m}$$
for all $f \in P_n$ and $g \in P_m$, and extended by linearity.  It is easy to see that this product corresponds to convolution of endomorphisms of the tensor algebra via the map~\eqref{E:block-duality}, when $\field=\C$.

For example,

\begin{picture}(410,80)(0,0)
\put(0,70){\circle*{3}}\put(0,55){\circle*{3}}\put(0,55){\line(1,0){15}}
\put(15,70){\circle*{3}}\put(15,55){\circle*{3}}\put(0,70){\line(1,0){15}}
\put(0,55){\line(0,0){15}}\put(15,55){\line(0,0){15}} \put(22,60){$\ast$}
\put(35,70){\circle*{3}}\put(35,55){\circle*{3}}\put(35,55){\line(1,1){15}}
\put(50,70){\circle*{3}}\put(50,55){\circle*{3}}\put(50,55){\line(-1,1){15}}
\put(60,60){$=$}
\put(75,70){\circle*{3}}\put(75,55){\circle*{3}}\put(75,55){\line(1,0){15}}
\put(90,70){\circle*{3}}\put(90,55){\circle*{3}}\put(75,55){\line(0,0){15}}
\put(105,70){\circle*{3}}\put(105,55){\circle*{3}}\put(75,70){\line(1,0){15}}
\put(120,70){\circle*{3}}\put(120,55){\circle*{3}}\put(90,55){\line(0,0){15}}
\put(105,55){\line(1,1){15}}\put(120,55){\line(-1,1){15}}
\put(127,60){$+$}
\put(145,70){\circle*{3}}\put(145,55){\circle*{3}}\put(145,70){\line(1,0){15}}
\put(160,70){\circle*{3}}\put(160,55){\circle*{3}}\put(145,55){\line(0,0){15}}
\put(175,70){\circle*{3}}\put(175,55){\circle*{3}}\qbezier(145,55)(160,45)(175,55)
\put(190,70){\circle*{3}}\put(190,55){\circle*{3}}\put(175,55){\line(-1,1){15}}
\put(160,55){\line(2,1){30}}\put(190,55){\line(-1,1){15}}
\put(197,60){$+$}
\put(215,70){\circle*{3}}\put(215,55){\circle*{3}}\put(215,70){\line(1,0){15}}
\put(230,70){\circle*{3}}\put(230,55){\circle*{3}}\put(215,55){\line(0,0){15}}
\put(245,70){\circle*{3}}\put(245,55){\circle*{3}}\qbezier(215,55)(237,45)(260,55)
\put(260,70){\circle*{3}}\put(260,55){\circle*{3}}\put(260,55){\line(-2,1){30}}
\put(230,55){\line(2,1){30}}\put(245,55){\line(0,0){15}}
\put(267,60){$+$}
\put(285,70){\circle*{3}}\put(285,55){\circle*{3}}\put(285,70){\line(1,0){15}}
\put(300,70){\circle*{3}}\put(300,55){\circle*{3}}\put(300,55){\line(1,0){15}}
\put(315,70){\circle*{3}}\put(315,55){\circle*{3}}\put(300,55){\line(-1,1){15}}
\put(330,70){\circle*{3}}\put(330,55){\circle*{3}}\put(315,55){\line(-1,1){15}}
\put(285,55){\line(3,1){45}}\put(330,55){\line(-1,1){15}}
\put(337,60){$+$}
\put(355,70){\circle*{3}}\put(355,55){\circle*{3}}\put(355,70){\line(1,0){15}}
\put(370,70){\circle*{3}}\put(370,55){\circle*{3}}\put(370,55){\line(-1,1){15}}
\put(385,70){\circle*{3}}\put(385,55){\circle*{3}}\qbezier(370,55)(385,45)(400,55)
\put(400,70){\circle*{3}}\put(400,55){\circle*{3}}\put(400,55){\line(-2,1){30}}
\put(355,55){\line(3,1){45}}\put(385,55){\line(0,0){15}}
\put(407,60){$+$}
\put(75,35){\circle*{3}}\put(75,20){\circle*{3}}\put(75,35){\line(1,0){15}}
\put(90,35){\circle*{3}}\put(90,20){\circle*{3}}\put(105,20){\line(1,0){15}}
\put(105,35){\circle*{3}}\put(105,20){\circle*{3}}\put(105,20){\line(-2,1){30}}
\put(120,35){\circle*{3}}\put(120,20){\circle*{3}}\put(120,20){\line(-2,1){30}}
\put(75,20){\line(3,1){45}}\put(90,20){\line(1,1){15}}
\end{picture}

\medskip

A \emph{breaking point} of a set partition $\calB$ is an integer $i\in \{0,1,\ldots, n\}$ for which there exists a subset $S\subseteq\calB$
such that
\[\bigcup_{B\in S} B =\{1,\ldots,i\} \text{ \ (and hence) \ } \bigcup_{B\in \calB\setminus S} B =\{i+1,\ldots,n\}\,.\]
Given a uniform block permutation  $f:\calA\to\calB$, let $B(f)$ denote the set of breaking points of $\calB$.
Note that $i=0$ and $i=n$ are  breaking points of any $f$. 
If $f$ is a permutation, that is if all blocks of $f$ are of size $1$, then  $B(f)=\{0,1,\ldots, n\}$. 

In terms of the diagram  of a uniform block permutation, if it is possible to put a vertical line between the first $i$
and the last $n-i$ vertices in the bottom row without intersecting an edge between the two sets of vertices, then $i$ is a breaking point. 

\begin{picture}(250,60)(0,0)
\put(30,25){$f=$}
\put(60,10){\red{\dashbox{1}(0,18){}}}
\put(70,40){\circle*{3}}\put(70,20){\circle*{3}}\put(70,20){\line(1,1){20}}
\put(80,10){\red{\dashbox{1}(0,18){}}}
\put(90,40){\circle*{3}}\put(90,20){\circle*{3}}\put(90,20){\line(3,1){60}}
\put(100,10){\red{\dashbox{1}(0,12){}}}
\put(110,40){\circle*{3}}\put(110,20){\circle*{3}}\put(110,20){\line(-2,1){40}}
\put(130,40){\circle*{3}}\put(130,20){\circle*{3}}\qbezier(110,20)(140,10)(170,20)
\put(150,40){\circle*{3}}\put(150,20){\circle*{3}}\put(130,20){\line(0,0){20}}
\put(170,40){\circle*{3}}\put(170,20){\circle*{3}}\put(150,20){\line(1,1){20}}
\put(180,10){\red{\dashbox{1}(0,18){}}}\put(190,40){\circle*{3}}\put(190,20){\circle*{3}}\put(170,20){\line(-3,1){60}}
\put(210,40){\circle*{3}}\put(210,20){\circle*{3}}\qbezier(70,40)(90,47)(110,40)
\put(190,20){\line(0,0){20}}\put(190,20){\line(1,0){20}}
\put(210,20){\line(0,0){20}}\put(190,40){\line(1,0){20}}
\put(220,10){\red{\dashbox{1}(0,18){}}}
\put(230,25){  $\Rightarrow \ \  B(f)=\{0,1,2,6,8\}$.}
\end{picture} 

\begin{lemma}\label{L:break}
If $i$ is a breaking point of $f$, then there exists a unique $(i,n-i)$-shuffle $\xi\in S_n$ and
unique uniform block permutations $f_{(i)} \in P_i$ and $f_{(n-i)}'\in  P_{n-i}$ such that
$$f=(f_{(i)}\times f'_{(n-i)})\cdot \xi^{-1}\,.$$
Conversely, if such a decomposition exists, $i$ is a breaking point of $f$.
\end{lemma} 
We illustrate this statement with an example where $i=4$ and 
$\xi=\mat{1 & 2 & 3 & 4 & 5 & 6\\2 & 3 & 5 & 6 & 1 & 4}$: 

\begin{center}
\begin{picture}(100,60)(0,-10)
\put(0,0){\circle*{3}}\put(0,40){\circle*{3}}\put(0,40){\line(4,-1){80}}
\put(20,0){\circle*{3}}\put(20,40){\circle*{3}}\put(20,40){\line(-1,-1){20}}
\put(40,0){\circle*{3}}\put(40,40){\circle*{3}}\put(40,40){\line(-1,-1){20}}
\put(60,0){\circle*{3}}\put(60,40){\circle*{3}}\put(60,40){\line(2,-1){40}}
\put(80,0){\circle*{3}}\put(80,40){\circle*{3}}\put(80,40){\line(-2,-1){40}}
\put(100,0){\circle*{3}}\put(100,40){\circle*{3}}\put(100,40){\line(-2,-1){40}}
\put(0,0){\dashbox{1}(0,20){}} \put(60,0){\dashbox{1}(0,20){}} \put(0,0){\dashbox{1}(60,0){}}\put(0,20){\dashbox{1}(60,0){}}
\put(80,0){\dashbox{1}(0,20){}} \put(100,0){\dashbox{1}(0,20){}} \put(80,0){\dashbox{1}(20,0){}}\put(80,20){\dashbox{1}(20,0){}}
\put(70,-10){\red{\dashbox{1}(0,20){}}}
\put(23,7){$f_{(4)}$} \put(83,7){$f'_{(2)}$}
\put(-35,17){$f=$} \put(110,25){$=\xi^{-1}$}
\end{picture}
\end{center}

\medskip

We are now ready to define the coproduct on $\calP$. Given $f\in  P_n$ set
$$\Delta(f):=\sum_{i\in B(f)}f_{(i)}\otimes f'_{(n-i)},$$
where $f_{(i)}$ and $f'_{(n-i)}$ are as in Lemma~\ref{L:break}. An example follows.
\begin{center}
\begin{picture}(320,80)(0,0)
\put(90,60){\large{$f = $}} 
\put(125,55){\circle*{3}}\put(125,70){\circle*{3}}\put(125,55){\line(1,1){15}}
\put(140,55){\circle*{3}}\put(140,70){\circle*{3}}\qbezier(125,55)(140,45)(155,55)
\put(155,55){\circle*{3}}\put(155,70){\circle*{3}}\put(140,55){\line(-1,1){15}}
\put(170,55){\circle*{3}}\put(170,70){\circle*{3}}\put(170,55){\line(-1,1){15}}
\put(185,55){\circle*{3}}\put(185,70){\circle*{3}}\put(170,55){\line(1,0){15}}
\put(200,55){\circle*{3}}\put(200,70){\circle*{3}}\put(185,55){\line(1,1){15}}
\put(200,55){\line(-1,1){15}}\qbezier(140,70)(155,77)(170,70)\qbezier(155,70)(177,77)(200,70)
\put(155,55){\line(1,1){15}}

\put(0,20){\large{$\Delta(f)=$}}
\put(50,20){\large{$f \otimes \emptyset \ + $}}
\put(105,15){\circle*{3}}\put(105,30){\circle*{3}}\put(105,15){\line(1,1){15}}
\put(120,15){\circle*{3}}\put(120,30){\circle*{3}}\qbezier(105,15)(120,5)(135,15)
\put(135,15){\circle*{3}}\put(135,30){\circle*{3}}\put(120,15){\line(-1,1){15}}
\put(135,15){\line(0,0){15}}\put(120,30){\line(1,0){15}}
\put(140,20){\large{$\otimes$}}
\put(155,15){\circle*{3}}\put(155,30){\circle*{3}}\put(155,15){\line(0,0){15}}
\put(170,15){\circle*{3}}\put(170,30){\circle*{3}}\put(155,15){\line(1,0){15}}
\put(185,15){\circle*{3}}\put(185,30){\circle*{3}}\put(170,15){\line(1,1){15}}
\put(185,15){\line(-1,1){15}}\qbezier(155,30)(170,37)(185,30)
\put(190,20){\large{$+$}}
\put(205,15){\circle*{3}}\put(205,30){\circle*{3}}\put(205,15){\line(1,1){15}}
\put(220,15){\circle*{3}}\put(220,30){\circle*{3}}\put(220,15){\line(-1,1){15}}
\put(235,15){\circle*{3}}\put(235,30){\circle*{3}}\qbezier(205,15)(220,5)(235,15)
\put(250,15){\circle*{3}}\put(250,30){\circle*{3}}\put(235,15){\line(1,1){15}}
\put(265,15){\circle*{3}}\put(265,30){\circle*{3}}\put(250,15){\line(-1,1){15}}
\put(250,15){\line(1,0){15}}\put(265,15){\line(0,0){15}}\qbezier(220,30)(235,37)(250,30)
\qbezier(235,30)(250,37)(265,30)
\put(270,20){\large{$\otimes$}}
\put(285,15){\circle*{3}}\put(285,30){\circle*{3}}\put(285,15){\line(0,0){15}}
\put(290,20){\large{$+$}}
\put(305,20){\large{$\emptyset\otimes f$.}}
\end{picture}
\end{center}
\vskip 0in  
Recall that an element $x\in\calP$ is called primitive if $\Delta(x)=x\otimes \emptyset +\emptyset\otimes x$. 
Every uniform block permutation with breaking set $\{ 0, n\}$ is primitive, but there other primitive elements in $\calP$. For example, the following element of $\field P_3$ is primitive:

 \begin{picture}(240,50)(0,0)
 \put(160,20){\circle*{3}}\put(160,35){\circle*{3}}\qbezier(160,35)(175,42)(190,35)
 \put(175,20){\circle*{3}}\put(175,35){\circle*{3}}\put(160,20){\line(0,0){15}}
 \put(190,20){\circle*{3}}\put(190,35){\circle*{3}}\put(160,20){\line(1,0){15}}
 \put(175,20){\line(1,1){15}}\put(190,20){\line(-1,1){15}}
 \put(200,25){$-$}
 \put(220,20){\circle*{3}}\put(220,35){\circle*{3}}\put(220,20){\line(1,0){15}}
 \put(235,20){\circle*{3}}\put(235,35){\circle*{3}}\put(220,20){\line(1,1){15}}
 \put(250,20){\circle*{3}}\put(250,35){\circle*{3}}\put(250,20){\line(-2,1){30}}
 \put(235,20){\line(1,1){15}}\put(235,35){\line(1,0){15}}
 \end{picture}

%

Recall that $\emptyset$ denotes the empty uniform block permutation. 
Let $\varepsilon:\calP\rightarrow \field$ be 
\[\varepsilon(f)=\begin{cases} 1 &\text{ if $f=\emptyset\in P_0$,}\\
0 & \text{ if $f\in  P_n$, $n\geq 1$. }
\end{cases}\]

\begin{theorem}
The graded vector space $\calP$, equipped with the product $\ast$, coproduct $\Delta$, unit $\emptyset$ and counit $\varepsilon$, is a  graded connected Hopf algebra. 
\end{theorem}
 Associativity  and coassociativity follow from basic properties of  shuffles (for the product one may also appeal to~\eqref{E:block-duality} and associativity of the convolution product).
  The existence of the antipode is guaranteed in any graded connected bialgebra.
Compatibility between $\Delta$ and $\ast$ requires a special argument. We sketch part of it. 

Let $\beta_{n,m}$ be the $(n,m)$-shuffle such that 
\[\beta_{n,m}(i)=\begin{cases} m+i & \text{  if $1\leq i\leq n$,}\\
i-n &  \text{  if $n+1\leq i\leq n+m$.}
\end{cases}\]
  The diagram of $\beta_{3,4}$ is shown below.

\begin{picture}(270,60)(0,10)
\put(150,60){\circle*{3}}\put(150,20){\circle*{3}}\put(150,20){\line(3,2){60}}
\put(170,60){\circle*{3}}\put(170,20){\circle*{3}}\put(170,20){\line(3,2){60}}
\put(190,60){\circle*{3}}\put(190,20){\circle*{3}}\put(190,20){\line(3,2){60}}
\put(210,60){\circle*{3}}\put(210,20){\circle*{3}}\put(210,20){\line(3,2){60}}
\put(230,60){\circle*{3}}\put(230,20){\circle*{3}}\put(230,20){\line(-2,1){80}}
\put(250,60){\circle*{3}}\put(250,20){\circle*{3}}\put(250,20){\line(-2,1){80}}
\put(270,60){\circle*{3}}\put(270,20){\circle*{3}}\put(270,20){\line(-2,1){80}}
\end{picture}

The inverse of $\beta_{n,m}$ is $\beta_{m,n}$.

Let $f\in P_n$, $g\in P_m$. A summand in $\Delta(f)\ast \Delta(g)$ is of the form \[\xi_1\cdot (f'\times g')\otimes \xi_2\cdot(f''\times g'')\]
where   $p\in B(f)$, $f'\in P_p$, $f''\in P_{n-p}$, $q\in B(g)$, $g'\in P_q$, $g''\in P_{m-q}$, $\xi_1\in \Sh(p,q)$, $\xi_2\in \Sh(n-p,m-q)$, 
 and there
 exist unique $\eta_1\in \Sh(p,n-p)$ and $\eta_2\in \Sh(q,m-q)$ such that $f\cdot \eta_1=f'\times f''$ and $g\cdot \eta_2=g'\times g''$.  
 
 Let $\beta:=1_p\times \beta_{n-p,q}\times 1_{m-q}$. Then
\begin{eqnarray*}
\xi_1\cdot (f'\times g')\times \xi_2\cdot (f''\times g'') &=& (\xi_1\times \xi_2)\cdot((f'\times g')\times (f''\times g''))\\ &  =& (\xi_1\times \xi_2)\cdot \beta \cdot((f'\times f'')\times(g'\times g''))\cdot \beta^{-1}\\
& =&   (\xi_1\times \xi_2)\cdot\beta\cdot (f\cdot \eta_1\times g\cdot\eta_2)\cdot \beta^{-1}\\
& =& (\xi_1\times \xi_2)\cdot\beta\cdot(f\times g)\cdot (\eta_1\times \eta_2)\cdot \beta^{-1}\,.
\end{eqnarray*}
Let $\xi:=(\xi_1\times \xi_2)\cdot \beta$ and $\eta:=(\eta_1\times \eta_2)\cdot \beta^{-1}$. One verifies that 
$\xi\in \Sh(n,m)$ and $\eta\in \Sh(p+q,n+m-p-q)$. Therefore, 
\[\xi_1\cdot (f'\times g')\times \xi_2\cdot (f''\times g'')=\xi\cdot (f\times g)\cdot \eta\]
is a summand in $\Delta(f\ast g)$.$\quad\Box$

\bigskip

Consider the following graded subspace of $\calP$: 
\[\calS:=\bigoplus_{n\geq 0}\field S_n\,.\]

\begin{prop} $\calS$ is a Hopf subalgebra of $\calP$. 
\end{prop}

$\calS$ is the  Hopf algebra of permutations of Malvenuto and Reutenauer~\cite{mr}.
Let $\sigma$ be a permutation. In the notation of~\cite{as}, the element $\sigma\in \calS$ corresponds to the
basis element $F_{\sigma}^*$ of $\SSym^*$, or equivalently the element 
$F_{\sigma^{-1}}$ of $\SSym$.

\subsection{Inverse monoid structure and self-duality}\label{S:selfdual}

As $\calS$, the Hopf algebra $\calP$ is self-dual. To see this, recall that a block permutation is a bijection $f:\calA\to\calB$ between two set partitions of $[n]$. Let $\Tilde{f}:\calB\to\calA$ denote the inverse bijection.
If $f$ is uniform then so is $\Tilde{f}$. The diagram of $\Tilde{f}\in P_n$ 
 is obtained by reflecting the diagram
of $f$ across a horizontal line.  Note that for $\sigma\in S_n\subseteq P_n$ we have $\Tilde{\sigma}=\sigma^{-1}$.

Let $\calP^*$ be the
graded dual space of $\calP$:
\[\calP^*=\bigoplus_{n\geq 0}(\field P_n)^*\,.\]
Let $\{f^*\mid f\in P_n\}$ be the basis of $(\field P_n)^*$ dual to the basis $P_n$ of $\field P_n$. 

\begin{prop} The map $\calP^*\to\calP$, $f^*\mapsto \Tilde{f}$, 
is an isomorphism of graded Hopf algebras.  
\end{prop}

The operation $f\mapsto\Tilde{f}$ is also relevant to the monoid
structure of $P_n$. Indeed, the following properties are satisfied
\[ f=f\Tilde{f}f \text{ \ and \ }\Tilde{f}=\Tilde{f}f\Tilde{f}\,.\]
Together with~\eqref{E:meet} below, these properties imply that $P_n$ is an {\em inverse monoid}~\cite[Theorem 1.17]{cp}.
The following properties are consequences of this fact~\cite[Lemma 1.18]{cp}:
\[\widetilde{fg}=\Tilde{g}\Tilde{f},\quad \Tilde{\Tilde{f}}=f\]
(they can also be verified directly).

\subsection{Factorizable monoid structure and the weak order}\label{S:weak}

Let $E_n$ denote the poset of set partitions of $[n]$: we say that $\calA\leq\calB$ if every bock of $\calB$ is contained in a block of $\calA$.
This poset is a lattice, and this structure is related to the monoid structure of uniform block 
permutations as follows. If $\id_{\calA}:\calA\to\calA$ denotes the uniform block permutation which is the identity map on the set of blocks of $\calA$, then
\begin{equation}\label{E:meet}
 \id_{\calA}\cdot\id_{\calB}=\id_{\calA\meet\calB}\,.
\end{equation}
In other words, viewing $E_n$ as a monoid under the meet operation $\meet$,
the map 
\[E_n\to P_n\,,\quad \calA\mapsto\id_{\calA}\,,\]
is a morphism of monoids. 

Any uniform block permutation $f\in P_n$ decomposes (non-uniquely) as \begin{equation}\label{E:factorization}
f=\sigma\cdot\id_{\calA} 
\end{equation}
for some $\sigma\in S_n$ and $\calA\in E_n$. Note that $\sigma$ is invertible and $\id_{\calA}$ is idempotent, by~\eqref{E:meet}. It follows that $P_n$ is a {\em factorizable inverse monoid}~\cite[Section 2]{ch}, ~\cite[Chapter 2.2]{l}. Moreover, by Lemma 2.1 in~\cite{ch}, any invertible element in $P_n$ belongs to $S_n$ and any idempotent element in $P_n$ belongs to (the image of) $E_n$. 
This lemma also guarantees that in~\eqref{E:factorization}, the idempotent $\id_{\calA}$ is uniquely determined by $f$ (which is clear since $\calA$ is
the domain of $f$). On the other hand, $\sigma$ is not unique, and we will make a suitable choice of this factor to define a partial order on $P_n$.

\begin{figure}[!h]
\centering
\includegraphics[width=10cm]{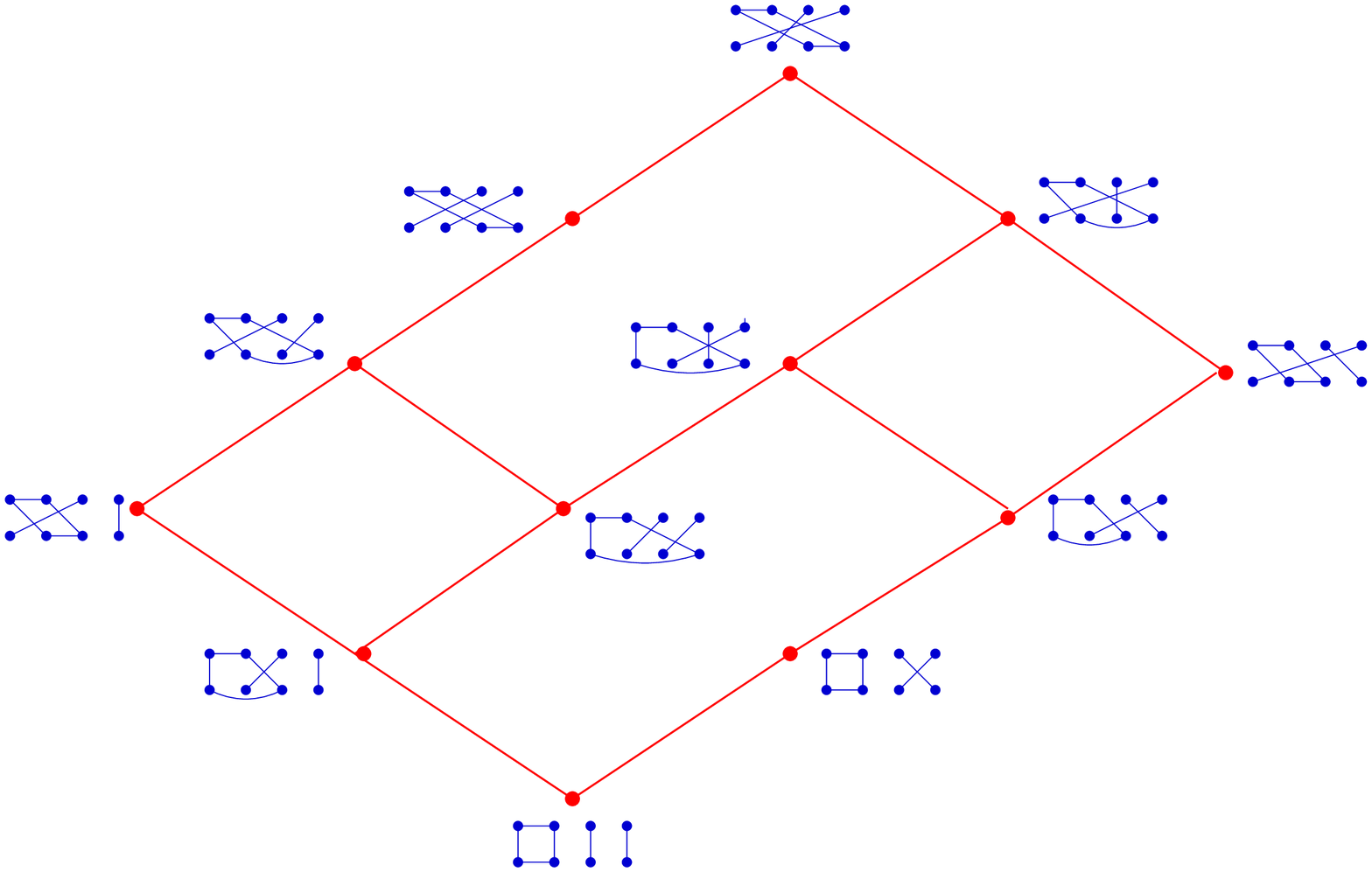}
\caption{The component of $P_4$ corresponding to $\calA=\{1,2\}\{3\}\{4\}$}\label{F:12-3-4} 
\end{figure}

\begin{figure}[!h]
\centering
\includegraphics[width=10cm]{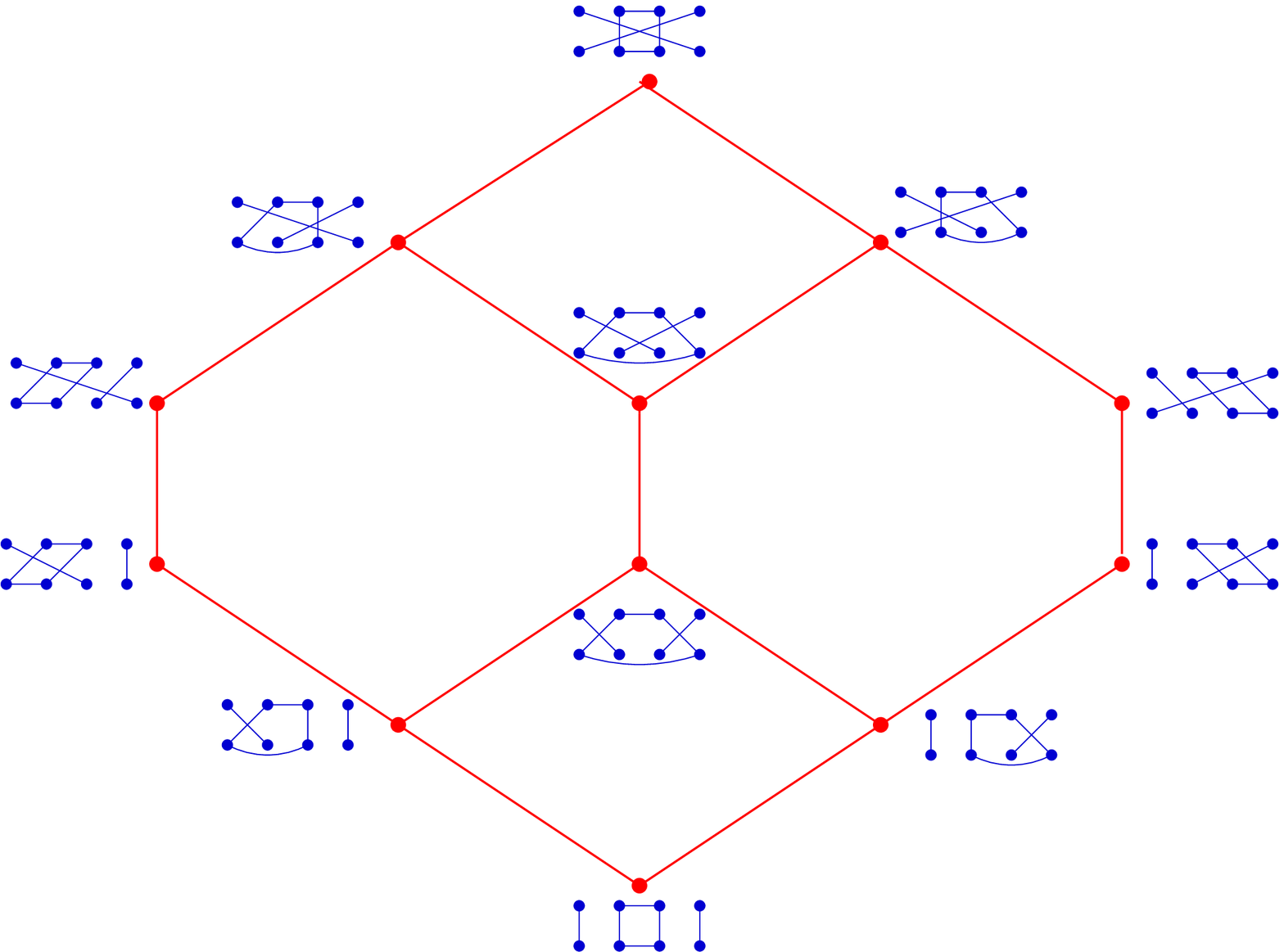}
\caption{The component of $P_4$ corresponding to $\calA=\{1\}\{2,3\}\{4\}$}\label{F:1-23-4} 
\end{figure}

\begin{figure}[!h]
\centering
\includegraphics[width=10cm]{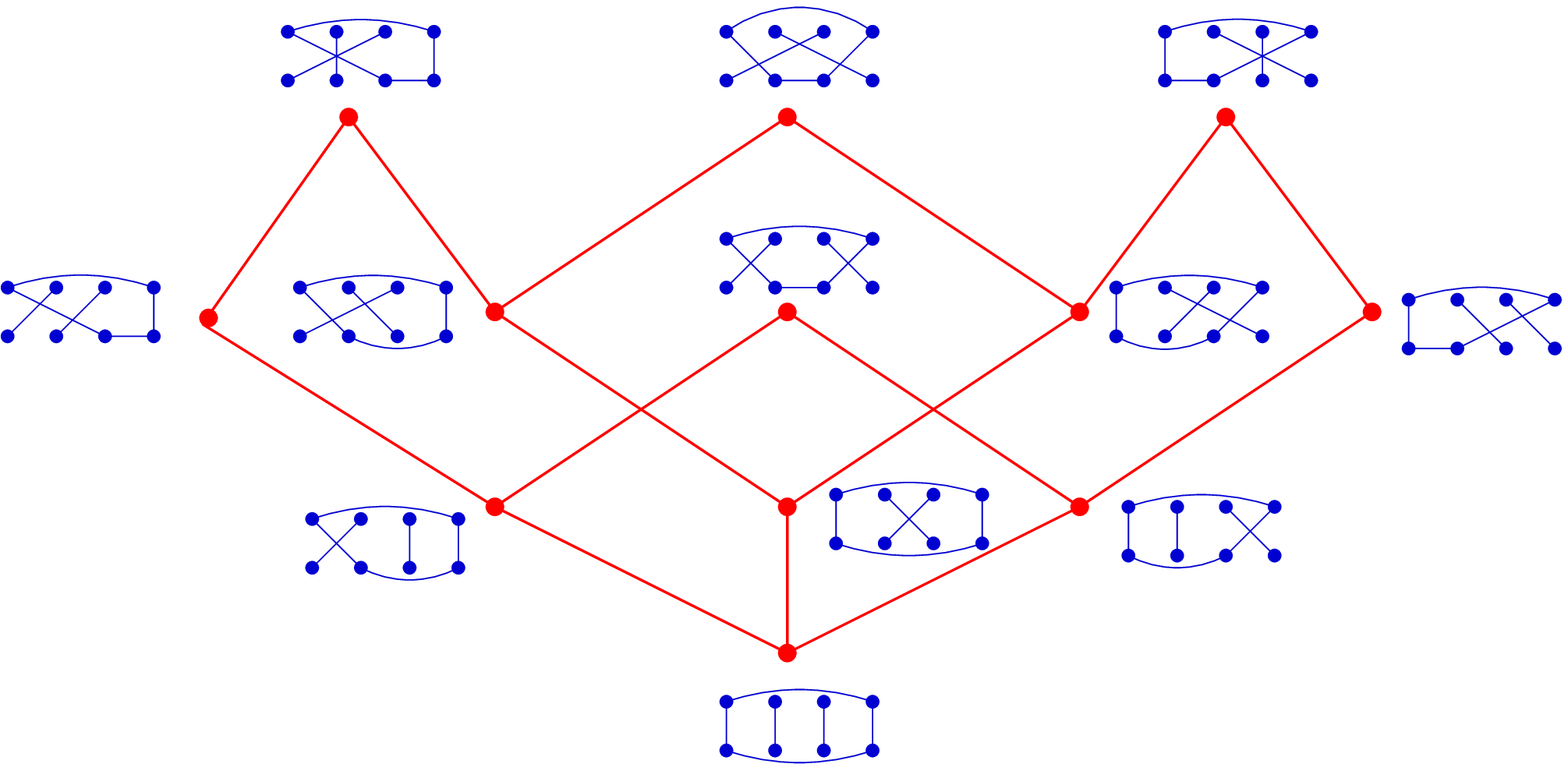}
\caption{The component of $P_4$ corresponding to $\calA=\{1,4\}\{2\}\{3\}$}\label{F:14-2-3} 
\end{figure}

\begin{figure}[!h]
\centering
\includegraphics[width=7cm]{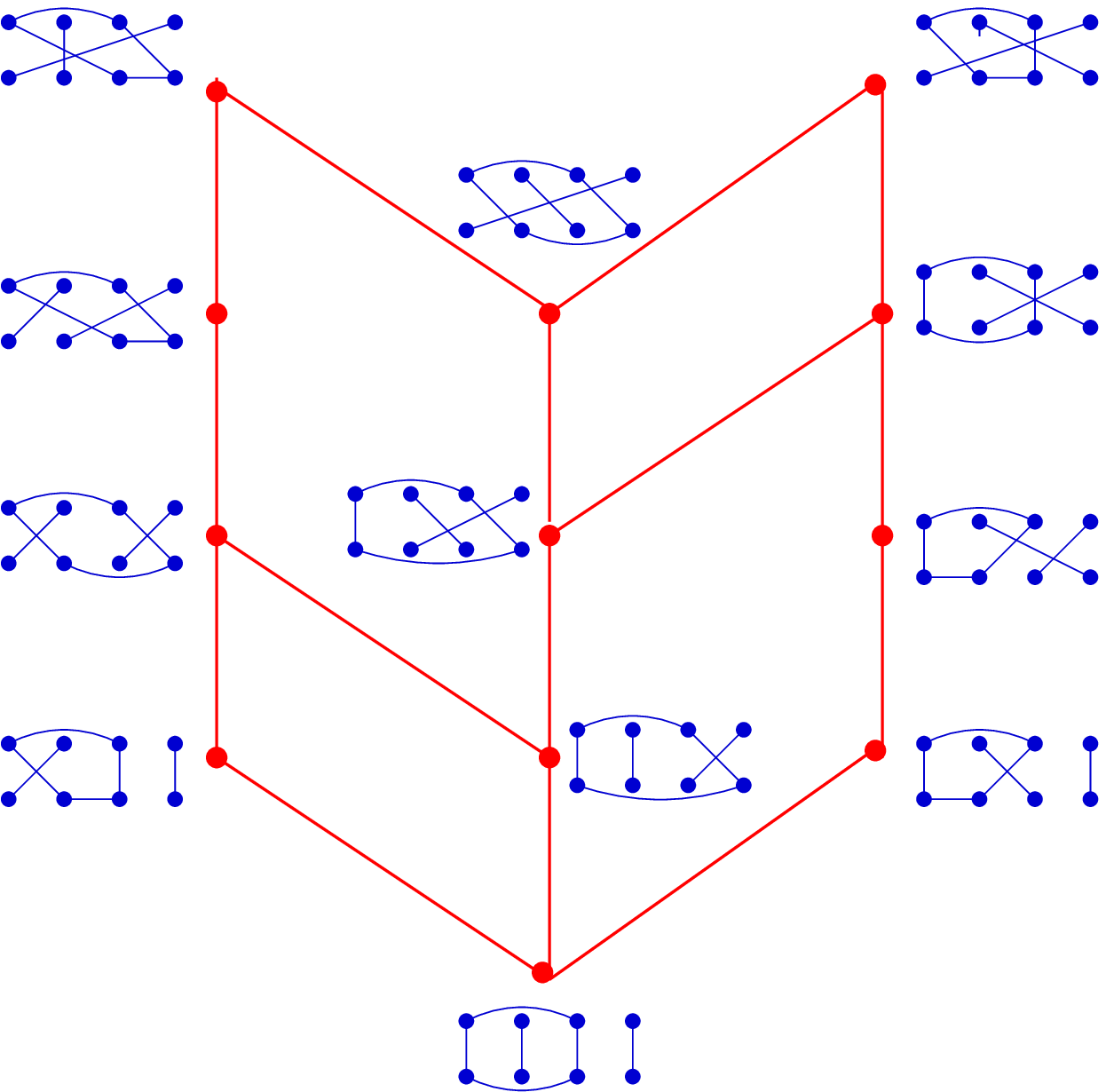}
\caption{The component of $P_4$ corresponding to $\calA=\{1,3\}\{2\}\{4\}$}\label{F:13-2-4} 
\end{figure}

\begin{figure}[!h]
\centering
\includegraphics[width=6cm]{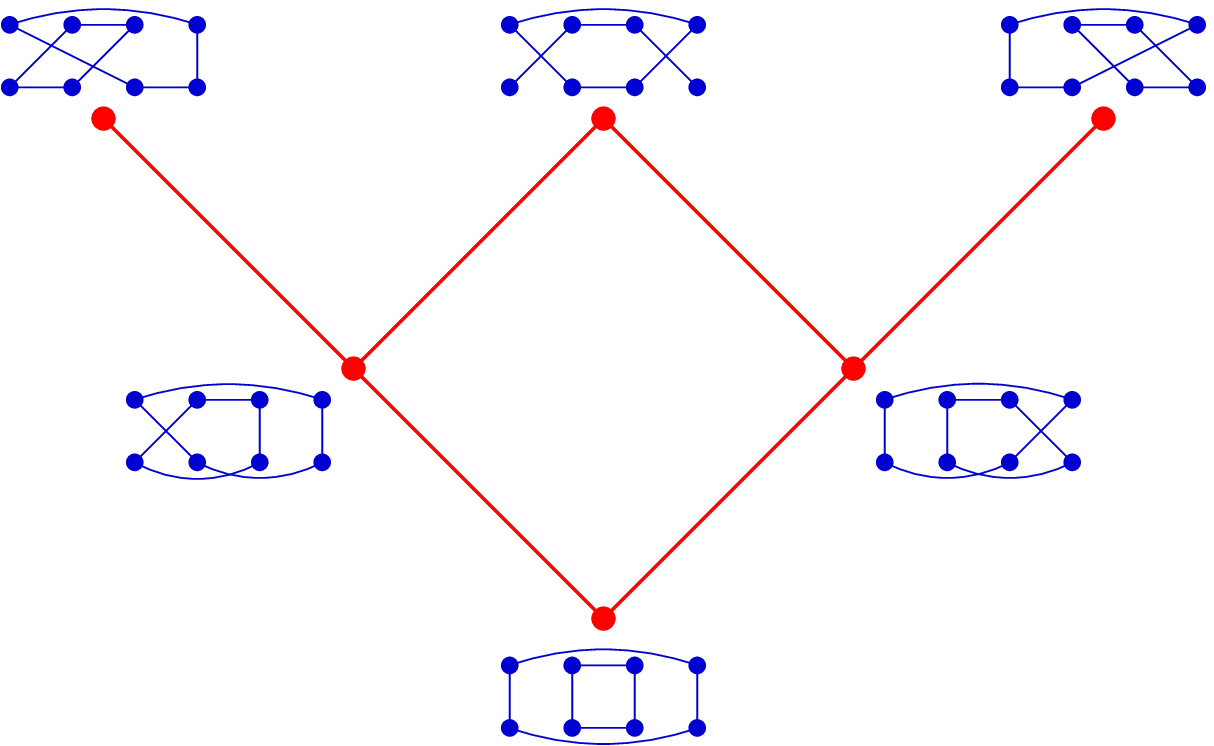}
\caption{The component of $P_4$ corresponding to $\calA=\{1,4\}\{2,3\}$}\label{F:14-23} 
\end{figure}

Consider the action of $S_n$ on $P_n$ by left multiplication. Given $\calA\in E_n$, the orbit of
$\id_{\calA}$ consists of all uniform block permutations $f:\calA\to\calB$ with
domain $\calA$, and the stabilizer is the {\em parabolic} subgroup
\[S_{\calA}:=\{\sigma\in S_n \mid \sigma(A)=A\ \forall A\in\calA\}\,.\]
Consider the set of $\calA$-\emph{shuffles}:
\[\Sh(\calA):=\{ \xi\in S_n \mid  \text{ if $i<j$ are in the same block of $\calA$ then $\xi(i)<\xi(j)$} \}\,.\]
It is well-known that these permutations form a set of representatives for the left cosets of the subgroup $S_{\calA}$. Therefore, given a uniform block permutation $f:\calA\to\calB$ there is a unique $\calA$-shuffle $\xi_f$ such that
\[f=\xi_f\cdot\id_{\calA}\,.\]
We use this decomposition to define a partial order on $P_n$ as follows:
\[f\leq g \iff \xi_f\leq \xi_g\,,\]
where the partial order on the right hand side is the left weak order on $S_n$ (see for instance~\cite{as}). We refer to this partial order as the {\em weak order} on $P_n$. Thus, $P_n$ is the disjoint union of certain subposets of the weak order on $S_n$:
\[P_n\cong \bigsqcup_{\calA\vdash[n]}\Sh(\calA)\]  
(in fact,  each $\Sh(\calA)$ is a lower order ideal $S_n$). Figures~\ref{F:12-3-4}-\ref{F:14-23} show 5 of the 15 components of $P_4$. Note that even when $\calA$ and $\calB$ are set partitions of the same type the posets $\Sh(\calA)$ and $\Sh(\calB)$ need not be isomorphic.

The partial order we have defined on $P_n$ should not be confused with the {\em natural partial order} which is defined on any inverse semigroup~\cite[Chapter 7.1]{cp2}, ~\cite[Chapter 1.4]{l}.

\vskip 0in

\vskip 0in 


\begin{rem} 
As observed by Sloane~\cite{sl}, there is a connection between
uniform block permutations and the {\em patience games} of Aldous and Diaconis~\cite{ad}. Starting from a deck of cards a patience game produces 
a number of card piles according to certain simple rules (the output is not unique). If the cards are numbered $1,\ldots,n$,
the initial deck is a permutation of $[n]$ and the resulting piles form a set partition of $[n]$. Suppose $\sigma\in S_n$. The set partitions $\calA$ such that
$\sigma\in\Sh(\calA)$ are precisely the possible outputs of patience games played from a deck of cards with $\sigma^{-1}(1)$ in the bottom, followed by $\sigma^{-1}(2)$, up to $\sigma^{-1}(n)$ on the top. Thus, 
uniform block permutations are in bijection with the pairs consisting of the input
and the output of a patience game via $(\sigma,\calA)\leftrightarrow \sigma\cdot\id_{\calA}$.
\end{rem}

 \subsection{The second basis and the Hopf algebra structure}\label{S:free}
 Following the ideas of~\cite{as}, we use the weak order on $P_n$ to define
 a new linear basis of the spaces $\field P_n$, on which the algebra
 structure of $\calP$ is simple.
 
For each element $g\in P_n$ let 
$$X_g:= \sum_{f\leq g} f\,.$$

By M\"obius inversion, the set $\{X_g \mid g\in P_n\}$ is a linear basis of $\calP_n$.

Given $p,q\geq 0$, let $\xi_{p,q}\in S_{p+q}$ be the permutation
\[\xi_{p,q}:=\begin{pmatrix} 1 & 2 & \ldots & p & p+1 & p+2&  \ldots &p+q\\
q+1 & q+2 & \ldots & q+p& 1 & 2 & \ldots & q
\end{pmatrix}\,.\]
This is the maximum element of $\Sh(p,q)$ under the weak order. The product of $\calP$ takes the following simple form on the $X$-basis.
\begin{prop} \label{P:X-product}
Let $g_1\in P_p$ and $g_2\in P_q$ be uniform block permutations. Then \[X_{g_1}\ast X_{g_2} = X_{\xi_{p,q}\cdot(g_1\times g_2)}\,.\] 
\end{prop} 

\begin{cor}\label{C:free} The Hopf algebra $\calP$ is free as an algebra and cofree as a graded coalgebra.
\end{cor}

Let $V$ denote the space of primitive elements of $\calP$. It follows that the
generating series of $\calP$ and $V$ are related by
\[\calP(x)=\frac{1}{1-V(x)}\,.\]
Since \[\calP(x)=1+x+3x^2+16x^3+131x^4+1496x^5+22482x^6+\cdots\]
we deduce that
\[V(x)=x+2x^2+11x^3+98x^4+1202x^5+19052x^6+\cdots\,.\]

\begin{rem} The same conclusion may be derived by introducing another basis
$$Z_g:= \sum_{f\geq g} f\,.$$
This has the property that
\[Z_{g_1}\ast Z_{g_2} = Z_{g_1\times g_2}\,.\] 
Note that $Z_{\id_{\calA}}$ is the element denoted $Z_{\calA}$ in Section~\ref{S:ideal}.
\end{rem}


\section{The Hopf algebra of symmetric functions in non-commuting variables}

Let $X$ be a countable set, the {\em alphabet}. A {\em word of length $n$} is a function $w:[n]\to X$. Let $\kX$ be the algebra of non-commutative
power series on the set of variables $X$. Its elements are infinite linear combinations of words, finitely many of each length, and the product is 
 concatenation of words.

The {\em kernel} of a word $w$ of length $n$ is the set partition $\calK(w)$ of $[n]$ whose blocks are the non-empty fibers of $w$.
Order the set of set partitions of $[n]$ by refinement, as in Section~\ref{S:weak}.  For each set partition $\calA$ of $[n]$, let
\[p_A:=\sum_{\calK(w)\leq A} w \in \kX\,.\]
This is the sum of all words $w$ such that if $i$ and $j$ are in the same block of $\calA$ then $w(i)=w(j)$. For instance
\[p_{\{1,3\}\{2,4\}}= xyxy+xzxz+ yxyx+ \cdots+ x^4+y^4+z^4+\cdots\,.\]

The subspace of $\kX$ linearly spanned by the  elements $p_\calA$ as $\calA$ runs over all set partitions of $[n]$, $n\geq 0$, is a subalgebra $\Pi$ of $\kX$, graded by length. The elements of $\Pi$ can be characterized  as those power series of finite degree that are invariant under any permutation of the variables.
$\Pi$ is the algebra of symmetric functions in non-commuting variables introduced by Wolf~\cite{w} and recently studied by Gebhard, Rosas, and Sagan~\cite{gs00, gs01,rs} in connection to Stanley's chromatic symmetric function.

$\Pi$ is in fact a graded Hopf algebra~\cite{brrz,am}. The coproduct is defined via evaluation of symmetric functions on two copies of the alphabet $X$. In order to describe
 the product and coproduct of $\Pi$ on the basis elements $p_\calA$ we introduce some notation.

Given set partitions  $\calA\vdash [n]$ and $\calB\vdash [m]$ let $\calA\times\calB$ be the set partition of $[n+m]$ whose blocks are the blocks of $\calA$ and the sets $\{b+n\, |\, b\in B\}$ where $B$ is a block of $\calB$.
This corresponds to the operation $\times$ on uniform block permutations in the sense that $\id_{\calA}\times\id_{\calB}=\id_{\calA\times\calB}$.
For example, if $\calA=\{1,3,4\}\{2,5\}\{6\}\vdash [6]$ and $\calB=\{1,4\}\{2\}\{3,5\}\vdash [5]$, then 
$\calA\times \calB =\{1,3,4\}\{2,5\}\{6\}\{7,10\}\{8\}\{9,11\}\vdash [11]$. 

To a set partition $\calA\vdash[n]$ and a subset $S\subseteq\calA$ we associate a new set partition $\calA_S$ as follows. Write
\[\bigcup_{A\in S} A=\{j_1,\cdots, j_m\} \subseteq [n]\]
 with $j_1<j_2<\cdots <j_m$. $\calA_S$ is the set partition of $[m]$ whose blocks
 are obtained from the blocks $A\in S$ by replacing each $j_t$ by $t$, for $1\leq t\leq m$. For instance, if $S=\{1,5\}\{2,7\}$ then $\calA_S=\{1,3\}\{2,4\}\vdash [4]$.
 
The product and coproduct of $\Pi$ are given by
\begin{align}
p_\calA p_\calB & = p_{\calA\times\calB}\,, \label{E:prodPi}\\
\Delta(p_\calA) & =\sum_{S\sqcup T=\calA}p_{\calA_S}\otimes p_{\calA_T}\,, \label{E:coprodPi}
\end{align}
the sum over all decompositions of $\calA$ into disjoint sets of blocks $S$ and $T$. For example, if $\calA=\{1,2,6\}\{3,5\}\{4\}$, then
\begin{eqnarray*}
\Delta(p_\calA) &=& p_\calA \otimes 1 + p_{\{1,2,5\}\{3,4\}}\otimes p_{\{1\}} 
+ p_{\{1,2,4\}\{3\}}\otimes p_{\{1,2\}} 
+p_{\{1,3\}\{2\}}\otimes p_{\{1,2,3\}} + \\  && p_{\{1,2,3\}}\otimes p_{\{1,3\}\{2\}} + 
p_{\{1,2\}} \otimes p_{\{1,2,4\}\{3\}} + p_{\{1\}}\otimes p_{\{1,2,5\}\{3,4\}} + 1 \otimes p_{\calA}\,.
\end{eqnarray*}

Consider now the direct sum of the subspaces $\calZ_n$ of $\field P_n$ introduced in Section~\ref{S:ideal}:
\[\calZ:=\bigoplus_{n\geq 0}\calZ_n \subset \calP \,.\]

\begin{theorem} $\calZ$ is a Hopf subalgebra of $\calP$. Moreover, the map
\[\calZ\to\Pi, \quad Z_{\calA}\mapsto p_{\calA}\]
is an isomorphism of graded Hopf algebras.
\end{theorem}

Thus the Hopf algebra of uniform block permutations $\calP$ contains the Hopf algebra $\Pi$ of symmetric functions in non-commuting variables. Note also that this reveals the
existence of a second operation on $\Pi$: according to Corollary~\ref{C:ideal}, each homogeneous component $\Pi_n$ carries an associative non-unital product that turns it into a right ideal of the
monoid algebra $\field P_n$. Connections between $\Pi$ and other combinatorial
Hopf algebras are studied in~\cite{am}.

\end{document}